\documentclass[conference]{IEEEtran}
\IEEEoverridecommandlockouts
\usepackage{cite}
\usepackage{amsmath,amssymb,amsfonts}
\usepackage{algorithmic}
\usepackage{graphicx}
\usepackage{hyperref}
\usepackage{textcomp}
\usepackage{xcolor}
\def\BibTeX{{\rm B\kern-.05em{\sc i\kern-.025em b}\kern-.08em
    T\kern-.1667em\lower.7ex\hbox{E}\kern-.125emX}}
\usepackage[normalem]{ulem}
\usepackage{multirow}
\usepackage{booktabs}
\usepackage{dblfloatfix}

\newcommand{\Q}{TPIA-Q}

\newcommand{\B}{TPIA-B}

\newcommand{\I}{TPIA-I}

\newcommand{\Isources}{slack current sources}

\newcommand{\phase}{\phi}

\newcommand{\squeezeup}{\vspace{-1.5mm}}
\newcommand{\squeezeuphalf}{\vspace{-.8mm}}

\begin{document}
\IEEEoverridecommandlockouts
\IEEEpubid{\makebox[\columnwidth]{\copyright 2023 IEEE \hfill}
\hspace{\columnsep}\makebox[\columnwidth]{ }}
\setlength{\abovedisplayskip}{3pt}
\setlength{\belowdisplayskip}{3.5pt}
\title{Actionable Three-Phase Infeasibility Optimization with Varying Slack Sources


\thanks{
* denotes the corresponding author}
}

\author{\IEEEauthorblockN{Elizabeth Foster$^{1*}$, Timothy McNamara$^1$, Amritanshu Pandey$^{2}$ and Larry Pileggi$^1$}
\IEEEauthorblockA{$^1$ Electrical \& Computer Engineering, \textit{Carnegie Mellon University}, Pittsburgh, PA} 
\IEEEauthorblockA{$^2$ Electrical \& Biomedical Engineering, \textit{University of Vermont}, Burlington, VT}
\{emfoster, tmcnama2, pileggi\}@andrew.cmu.edu, {amritanshu.pandey}@uvm.edu} 

\maketitle

\begin{abstract}
Modern distribution grids that include numerous distributed energy resources (DERs) and battery electric vehicles (BEVs) will require simulation and optimization methods that can capture behavior under infeasible operating scenarios to assess reliability. A three-phase infeasibility analysis (TPIA) localizes and identifies power deficient areas in distribution feeders via a non-convex optimization that injects and subsequently minimizes slack sources, subject to AC network constraints.  In this paper, we extend the TPIA framework by introducing operational bounds to ensure realistic, actionable solutions. We incorporate current, reactive power, and susceptance slack sources to model real-world assets, and discuss their potential use cases.  
 We show that the voltage-bounded TPIA formulations provide actionable solutions for realistic networks of up to 5360 nodes where power flow simulations either fail or return low-voltage solutions. We demonstrate reactive power compensation using the slack susceptance formulation on an infeasible test case.
\end{abstract}


\section{Introduction}
Modern electric grids are seeing an influx of distributed energy resources (DERs) and smart technologies like battery electric vehicles (BEVs) at the distribution level \cite{IEA} \cite{NAS} \cite{derreport}. These new resources and technologies present opportunities to achieve climate goals, amongst other benefits; however, they also present new challenges as system planners attempt to build out the infrastructure required to incorporate these assets while maintaining grid reliability and resiliency \cite{NAS} \cite{derreport}. 

System planners use tools generally limited to traditional methods like three-phase power flow, which outputs nodal voltages and line/transformer flows under different feasible operating conditions \cite{derreport}. Three-phase power flow will struggle to adequately assess the challenges presented by BEV and DER growth because scenarios requiring upgrades will fail to converge for infeasible networks that cannot satisfy AC network constraints. For example, simulation alone cannot identify optimal locations for necessary infrastructure upgrades due to rapid electrification or lack of reactive power in weak areas of the network. Currently available tools for distribution grid power flow offer limited, if any, optimization capabilities \cite{derreport}. New distribution planning tools and capabilities are needed that are capable of running optimizations while satisfying three-phase AC unbalanced network constraints. These new tools will help system planners make optimal decisions on updating and installing new infrastructure when combined with operational, financial, and policy constraints \cite{SCE}.

While the availability of such complete three-phase power flow tools is limited, three-phase optimization is an active and broad research area \cite{DGcommittee}. Several optimization approaches are targeted to three-phase optimal power flow (OPF), which typically minimizes economic dispatch subject to network constraints. These OPF methods have often been formulated using power balance network equations as seen in \cite{gill}, which applies OPF to active distribution networks. The authors in \cite{Dall} extend traditional power-balance constrained OPF to include slack variables that indicate constraint violations when there is no feasible solution. Other three-phase optimizations include technical objectives, like those that minimize losses \cite{Ochoa}, while satisfying network equality constraints. 

Formulations using power-balance equations \cite{gill} \cite{Dall} \cite{Ochoa} are highly nonlinear and therefore do not scale well in general. This behavior is shown in \cite{Claeys}, which developed a four-wire OPF with both current-voltage (I-V) and power balance formulations and found that power-balance did not converge for 35\% of test cases. Not only was I-V faster across the moderately sized test cases where power balance did converge, but the I-V formulation also returned more practical solutions. In the power balance formulation for 4-wire systems, some of the nodes were non-physically grounded, thereby resulting in erroneous voltage solutions \cite{Claeys}. Furthermore, many of these three-phase network-constrained optimizations will fail when there is no feasible power flow solution if slack variables or constraint violations are not considered.

A three-phase infeasibility analysis (TPIA) is another type of three-phase optimization formulated subject to current-based network constraints for use in distribution grid planning. TPIA localizes and identifies deficient areas, in addition to the \textit{amount} of the deficiency, in situations where power flow fails \cite{PSCC}. At a high level, it is a non-convex optimization study that identifies scenarios where AC three-phase power flow constraints are infeasible by using a set of variable slack sources with a physical meaning to compensate for \textit{missing power} in the system. In TPIA, solutions with nonzero slack sources occur when power flow would otherwise fail. 

Initially, the concept of quantifying infeasibility used power balance network constraints for positive sequence systems to find a solvable boundary for otherwise insolvable power flow simulations \cite{overbye}. We formulated an \textbf{infeasibility analysis} using \Isources{} and KCL network constraints for transmission networks in \cite{marko} and three-phase networks in \cite{PSCC}, where we additionally developed an L1-norm objective formulation that effectively localized the sources of infeasibility. Our work in \cite{PSCC} and \cite{marko} focused only on slack current models that can be directly added into the KCL network constraints. These works, \cite{PSCC} and \cite{marko}, did not include analysis on the selection of \Isources{} or other realistic physical asset models, such as slack reactive power or slack susceptance sources, that naturally represent physical assets.  In addition, operational constraints, like nodal voltages bounds, were not included in the prior formulation. 

In this paper, we develop a TPIA formulation in Section \ref{methodology} with realistic, actionable solutions by incorporating bus voltage magnitude restrictions using inequality constraints and introducing slack sources that model real-world devices. 
We consider the relationship between physical assets that address real-world challenges, like capacitor banks, and subsequently develop TPIA formulations for reactive power compensation using slack reactive power or susceptance sources. 
We empirically compare TPIA with standard three-phase power flow simulations and demonstrate implementing reactive power compensation based on a TPIA solution to make a formerly-unsolvable test case feasible in Section \ref{result}.

\section{Background} \label{background}
In this section, we discuss the infeasibility analysis concept and the circuit-simulation heuristics necessary for robust convergence of large-scale three-phase optimization.


\subsection{Infeasibility Analyses} \label{infeasbackground} 
Infeasibility analysis was first introduced for transmission networks in \cite{overbye}, and later in terms of an equivalent circuit for power grids in \cite{marko}. \cite{marko} introduced \Isources{} at each node in the positive sequence circuit model of the power grid. It minimized the norm of the \Isources{} to identify when positive sequence power flow was infeasible and returned the amount of missing current at each node in the system. \cite{Cindy} extended \cite{marko} to localize the sources of infeasibility by using an L1 regularization term in a least-squares objective. In \cite{PSCC}, we introduced a scalable TPIA framework using three-phase current-based network constraints with two objective functions: one using a least-squares minimization; and the other using an L1-norm, which localized the solution to a subset of nodes.

However, the equivalent circuit approaches in \cite{PSCC} \cite{marko} \cite{Cindy} do not relate slack source values with realizable physical assets on the grid. Furthermore, they do not enforce physical bounds on the grid voltages, thereby producing results with non-actionable infeasible solutions. Additionally, minimizing current produces solutions where the net power injection might be impractical given the varying base voltages across the distribution grid (e.g. providing a 1 amp current injection at a 13.2 kV bus requires much more power than supplying 1 amp at a 240 V bus). In this paper, we address these limitations by: (a) formulating TPIA with alternative slack sources (slack susceptance and reactive power); (b) restructuring the objective function to minimize net complex power for all slack formulations; and (c) enforcing operational bounds using inequality constraints. 
\squeezeup
\subsection{Circuit Heuristics for Three-phase Optimization} \label{heuristics}
To ensure robust convergence of TPIA when there are inequality constraints from operational limits, we leverage a limiting technique for Newton-Raphson (NR) that is based on our understanding of individual device model characteristics. First developed for diodes, this heuristic dampens the update steps of primal and dual variables associated with inequality constraints so that their respective values never enter infeasible regions \cite{diode}. 

\squeezeup
\squeezeup
\section{Methodology} \label{methodology}

In TPIA, if the network is feasible, the slack sources have values of zero at convergence. Otherwise, nonzero slack sources occur and can be used to identify and localize grid infeasibility for planning insight or corrections during operation. To ensure that various TPIA formulations provide meaningful solutions, we introduce operational limits in Section \ref{limit}. We consider various slack formulations in order to discuss their practical applicability in grid operation and planning. The TPIA formulations with different slack sources are: \I{} with \Isources{}; \Q{} with slack reactive power sources; and \B{} with slack susceptance sources. All of these formulations provide flexibility to impose upper and lower limits on the slack sources (e.g., constraining a slack susceptance $B$ per phase $\phase$ by its upper rating: $B^\phase_s\leq\overline{B^\phase_s}$ \cite{shunt_max}). Our prior work in \cite{PSCC} covered a version of \I{}, but  we reproduce it here for completeness.
\squeezeuphalf
\subsection{Three-Phase Power Flow Constraints}
We use nonlinear KCL network equations to represent the steady-state grid physics in the TPIA formulations. Due to the modular circuit-based approach, the TPIA framework can model any three-wire or four-wire distribution feeders with various transformer configurations and distribution assets like inverter-based distributed generators, voltage regulators, triplex loads, capacitors, fuses, shunts, and switches. The models for these devices are detailed in \cite{Amrit}. 

We describe the KCL-based three-phase AC network constraints in \eqref{eqn: PowerFlow-3} and \eqref{eqn: PowerFlow-4}, where $f^\phase_r(x)$ and $f^\phase_i(x)$, are the respective rectangular real and imaginary nodal current equations. $x$ is a vector of power flow variables that includes the real and imaginary voltages, $V^\phase_R$ and $V^\phase_I$, per phase $\phase$. $G$ and $B$ are matrices of network conductances and susceptances. $I^{NL}_R$ and $I^{NL}_I$ are nonlinear currents from either PV or PQ buses, but for brevity, we represent them as though they are from PQ buses and omit the reactive power variable and voltage set-point equation for PV buses. Given a solution $x$ for a feasible network, $f^\phase_r(x)$ and $f^\phase_i(x)$ are equal to 0. For infeasible networks, values for $f^\phase_r(x)$ and $f^\phase_i(x)$  are non-zero.
\begin{align}
    f_r^\phase(x) = G V^{\phase}_R - B V^{\phase}_I + I^{NL}_R(V^\phase_R, &V^\phase_I)  = 0
    \label{eqn: PowerFlow-3} \\
    f_i^\phase(x) = B V^{\phase}_R + G V^{\phase}_I + I^{NL}_I(V^\phase_R, &V^\phase_I) = 0
    \label{eqn: PowerFlow-4} 
\end{align}

The nonlinear currents $I^{NL}_R$ (\ref{eqn: PowerFlow-1}) and $I^{NL}_I$ (\ref{eqn: PowerFlow-2}) are vectors and functions of power injections. We later define slack reactive power sources using modified versions of (\ref{eqn: PowerFlow-1}) and (\ref{eqn: PowerFlow-2}).
\begin{align}
I^{NL}_R(V^{\phase}_R, V^{\phase}_I) = (P^{\phase}\odot &V_R^{\phase}  + Q^{\phase}\odot V_I^{\phase}) \nonumber\\ 
   & \oslash ((V_R^{\phase})^2 + (V_I^{\phase})^2)  \label{eqn: PowerFlow-1} \\
    I^{NL}_I(V^{\phase}_R, V^{\phase}_I) = (P^{\phase}\odot &V_I^{\phase}  - Q^{\phase}\odot V_R^{\phase}) \nonumber \\
   & \oslash ((V_R^{\phase})^2 + (V_I^{\phase})^2)   \label{eqn: PowerFlow-2}
\end{align}
\subsection{Operational Limits} 
\label{limit}

Prior work has focused on how slack sources can be used to satisfy the equality constraints that represent the physics of a power flow problem (i.e. the power or current balance at each node). However, an infeasibility analysis that is only constrained by the AC network equations may produce a solution that satisfies KCL but contains voltage magnitudes or line flows outside of acceptable ranges. This solution would thereby render any information non-actionable.

To address this limitation, inequality constraints to enforce operational limits can be imposed by the user in conjunction with infeasibility sources to find practical solutions. In this work, we impose a constant-valued upper bound ($\overline{V_i}$)  and lower bound ($\underline{V_i}$) on the complex voltage magnitudes for each phase $\phase\in\{A,B,C\}$ of each bus $i \in n$:
\begin{subequations}
\begin{align}
\overline{l}^\phase_i(x): ~(V_{R,i}^\phase)^2 + (V_{I,i}^\phase)^2 - (\overline{V_i})^2 \leq 0 \label{eqn:vlimupper} \\
\underline{l}^\phase_i(x): ~(V_{R,i}^\phase)^2 + (V_{I,i}^\phase)^2 -  (\underline{V_i})^2 \geq 0 
\label{eqn:vlimlower}
\end{align}  
\end{subequations} 

$\overline{l}^\phase_i(x)$ and $\underline{l}^\phase_i(x)$ represent the voltage bounds as general inequality functions and will be used in TPIA definition given in Section \ref{generalTPIA}. Additional inequality constraints enforcing upper bounds on current across branches or power through transformers can be defined and implemented similarly.
\subsection{General TPIA Format}
\label{generalTPIA}
The general format for TPIA is the same for all slack variable formulations and is described in (\ref{eqn: gen}). 
\begin{subequations}
\label{eqn: gen}
\begin{align}
\min \quad & \sum_{m=1}^n \sum_{\phase \in \{A,B,C\}} \alpha g^{\phase}_m(s) \label{eqn: genvcobj} \\
\textrm{s.t.} \quad & f_r^{\phase}(x) - h_r^{\phase}(s)= 0 
\label{eqn: genvcrealeqn}\\
& f_i^\phase(x) - h_i^{\phase}(s) = 0 \label{eqn: genvcimageqn}
\\
& 
\overline{l}^\phase_i(x) \leq 0  \label{eqn:geninequpper} \\
&
\underline{l}^\phase_i(x) \geq 0 \label{eqn:genineqlower}
\end{align}
\end{subequations}
The objective function in \ref{eqn: genvcobj} is the sum of squares of the complex power injection, $g_m^{\phase}(s)$, from the set of slack variables $s$ at all $n$ buses and all connected phases $\phase$. We select this objective function to address some of the challenges discussed in Section \ref{infeasbackground}. 
The scaling parameter $\alpha$ is tuned to bring the magnitudes of objective function gradients into balance with magnitudes of inequality constraint gradients; tuning $\alpha$ properly reduced oscillatory behavior in the solver and improved convergence speed.
Each formulation is subject to KCL network constraints modified by $h^\phase_r(s)$ and $h^\phase_i(s)$, which represent real (\ref{eqn: genvcrealeqn}) and imaginary (\ref{eqn: genvcimageqn}) current injections from slack sources. The functions $g_m^{\phase}(s), h^{\phase}_r(s)$, and $h^{\phase}_i(s)$ will change depending on the choice of slack source; the voltage bounds in (\ref{eqn:geninequpper}) and (\ref{eqn:genineqlower}) are the same for all formulations. 
 
We now develop \I{} using slack current sources.
The function $g^{\phase}_m(s)$ in (\ref{eqn: currobj}) represents the square of complex power using voltage and current variables. 
This objective ensures that the formulation targets missing power at a given node and is voltage-invariant, unlike our prior work in \cite{PSCC}. The current injections in (\ref{eqn: currreal}) and (\ref{eqn: currimag}) are simply the \Isources{}. 
\squeezeup
\squeezeup
\begin{subequations} \label{eqn:TPIA_I}
\begin{align}
g_m^{\phase}(s) &= \left((i_{fR,m}^{\phase})^2 +  (i_{fI,m}^{\phase})^2\right) \left((V_{R,m}^{\phase})^2 +  (V_{I,m}^{\phase})^2\right)\label{eqn: currobj} \\
h_{r}^\phase(x) &= i_{fR}^{\phase} \label{eqn: currreal} \\
h_{i}^\phase(x) &= i_{fI}^{\phase} \label{eqn: currimag}
\end{align}
\end{subequations}
We next address some limitations of \I{} by using alternative slack sources to formulate TPIA. 

\squeezeuphalf
\begin{table*}[ht!]
\centering
\caption{\label{table:pq} Comparison of infeasibility solver to standard powerflow}
\begin{tabular}{|cc|ccc|ccccc|}
\hline
\multicolumn{2}{|c|}{}                                   & \multicolumn{3}{c|}{\textbf{Three-phase Power Flow}}                                                           & \multicolumn{5}{c|}{\textbf{Voltage-Bounded Three-phase Infeasibility Analysis}}                                                                                \\ \hline
\multicolumn{1}{|c|}{\textbf{Test case}}  & \textbf{\# Nodes} & \multicolumn{1}{c|}{\textbf{\# $V_{mag}$ Viol.}} & \multicolumn{1}{c|}{\textbf{Iter.}} & \textbf{Time (s)} & \multicolumn{1}{c|}{\textbf{ \# $V_{mag}$ Viol.}} & \multicolumn{1}{c|}{\textbf{ $\Sigma |P_{feas}|$ (kW)}} & \multicolumn{1}{c|}{\textbf{ $\Sigma |Q_{feas}|$ (kVAR)}} & \multicolumn{1}{c|}{\textbf{Iter.}} & \textbf{Time (s)} \\ \hline
\multicolumn{1}{|c|}{R1-12.47-3\_OV} & 76                & \multicolumn{1}{c|}{\textcolor{red}{DNC}}                             & \multicolumn{1}{c|}{\textcolor{red}{DNC}}            & \textcolor{red}{DNC}               & \multicolumn{1}{c|}{0} & \multicolumn{1}{c|}{413.84}              & \multicolumn{1}{c|}{534.83}                & \multicolumn{1}{c|}{9}              & 0.47              \\ \hline
\multicolumn{1}{|c|}{R2-25.00-1\_OV} & 800               & \multicolumn{1}{c|}{\textcolor{red}{DNC} }                             & \multicolumn{1}{c|}{\textcolor{red}{DNC}}            & \textcolor{red}{DNC}               & \multicolumn{1}{c|}{0} & \multicolumn{1}{c|}{8,718.91}            & \multicolumn{1}{c|}{16,308.32}             & \multicolumn{1}{c|}{26}             & 12.42             \\ \hline
\multicolumn{1}{|c|}{R4-12.47-1}     & 1599              & \multicolumn{1}{c|}{16}                              & \multicolumn{1}{c|}{4}              & 0.618772          & \multicolumn{1}{c|}{0} & \multicolumn{1}{c|}{665.69}              & \multicolumn{1}{c|}{1,168.13}              & \multicolumn{1}{c|}{35}             & 32.17             \\ \hline
\multicolumn{1}{|c|}{R1-12.47-1}     & 1833              & \multicolumn{1}{c|}{135}                             & \multicolumn{1}{c|}{4}              & 0.69          & \multicolumn{1}{c|}{0} & \multicolumn{1}{c|}{631.82}              & \multicolumn{1}{c|}{1,202.72}              & \multicolumn{1}{c|}{20}             & 21.66             \\ \hline
\multicolumn{1}{|c|}{R5-12.47-3}     & 4046              & \multicolumn{1}{c|}{2912}                            & \multicolumn{1}{c|}{4}              & 1.57          & \multicolumn{1}{c|}{0} & \multicolumn{1}{c|}{2,311.98}            & \multicolumn{1}{c|}{4,854.95}              & \multicolumn{1}{c|}{30}             & 70.68             \\ \hline
\multicolumn{1}{|c|}{R3-12.47-3}     & 5360              & \multicolumn{1}{c|}{349}                             & \multicolumn{1}{c|}{4}              & 2.53          & \multicolumn{1}{c|}{0} & \multicolumn{1}{c|}{589.61}              & \multicolumn{1}{c|}{622.66}                & \multicolumn{1}{c|}{32}             & 124.89            \\ \hline
\end{tabular} 
\squeezeup
\squeezeup
\squeezeuphalf
\end{table*}
\subsection{TPIA for Reactive Power Compensation}
System planners might want to explore installing either static synchronous compensators (STATCOMs) or capacitor banks for reactive power compensation. However, \I{} cannot restrict injections to reactive power alone without imposing additional equality constraints that limit real power. Therefore we instead expand TPIA to use slack reactive power sources (\Q) or slack susceptance sources (\B). Both of these formulations model real-world reactive power compensation devices without adding new constraints to TPIA.

The equations in \ref{eqn:PQ} describe \Q{}, which evaluates infeasibility using slack reactive power sources instead of \Isources{}. The function $g_m^\phase(s)$ in (\ref{eqn: PQ_obj}) represents the square of the slack reactive power, $Q_s$, at bus $m$ and phase $\phase$. The KCL current injections from the slack sources, $I^s_R$ (\ref{eqn: PQreal}) and $I^s_I$ (\ref{eqn: PQimag}), are equivalent to (\ref{eqn: PowerFlow-1}) and (\ref{eqn: PowerFlow-2}) except the constant terms $P$ and $Q$ from those equations are replaced with a constant 0 and variable $Q_s$, respectively.
\begin{subequations}
\label{eqn:PQ}
\begin{align}
g_m^{\phase}(s) &= (Q_{s,m}^{\phase})^2 \label{eqn: PQ_obj} \\
h_{r}^{\phase}(x) &= I^{s}_{R}(V^{\phase}_{R}, V^{\phase}_{I}, 0, Q_{s}^{\phase}) \label{eqn: PQreal}\\
h_{i}^{\phase}(x) &= I^{s}_{I}(V^{\phase}_{R}, V^{\phase}_{I}, 0, Q_{s}^{\phase}) \label{eqn: PQimag} 
\end{align}
\end{subequations}


Depending on the type of asset a planner is considering installing, slack susceptance sources could also represent reactive power compensation in TPIA.
This formulation is detailed in (\ref{eqn: GB}) and is referred to as \B{}.
The function $g_m^\phase(s)$ in (\ref{eqn: Bobj}) is still the square of reactive power at bus $m$ and phase $\phase$, but it is now formulated with bus voltages and slack susceptance $B_s$. 
The current injections from slack susceptance are functions of voltage and susceptance and are shown in (\ref{eqn: Breal}) and (\ref{eqn: Bimag}).
\begin{subequations}
\label{eqn: GB}
\begin{align}
g_m^{\phase}(s) &= (B_{s,m}^{\phase})^2\left( (V^{\phase}_{R,m})^2 + (V^{\phase}_{I,m})^2 \right)^2 \label{eqn: Bobj} \\
h_r^{\phase}(x) & = B_s^{\phase} \odot V^{\phase}_I   
\label{eqn: Breal}\\
h_i^{\phase}(x) &= - B_s^{\phase} \odot V^{\phase}_R  \label{eqn: Bimag} 
\end{align}
\end{subequations}

When minimizing the sum of squares of complex power, \Q{} and \B{} become mathematically equivalent.
Adjusting the objective function to represent different real-world needs (e.g. minimizing the amount of slack \textit{capacitance} installed) could lend itself to unique use cases and solutions for system planners to exploit.
Furthermore, if a planner wanted to directly model real power slack sources in TPIA, the modifications to (\ref{eqn:PQ}) to include a $P^{\phi}_{s,m}$ variable are straightforward. As an example, this formulation could be used when determining where to site grid-connected batteries.

\subsection{Optimization Solution Quality}
For each of the above formulations, we solve the optimization problem using the primal-dual interior point (PDIP) approach \cite{Boyd} aided by the circuit-simulation heuristics from Section \ref{heuristics}. 
We seek a local minimizer by solving for a set of primal and dual variables that satisfy the perturbed first-order KKT conditions of the Lagrangian. 
The underlying problem is non-convex due to nonlinear network constraints, so no standard approach exists to obtain the global optima. 
However, our prior works in \cite{Amrit} and \cite{two_stage} have proposed heuristics that can be used to improve convergence to practical and meaningful solutions and to avoid saddle points.

\section{Experimental Set-up} \label{experiment}
To evaluate TPIA, we incorporated the new formulations into our three-phase power grid implementation used in \cite{PSCC} in order to use circuit heuristics from Section \ref{background} \cite{code}.
The implementation is written in Python 3, and all experiments were run on a 2019 MacBookPro with a 2.6 GHz 6-Core Intel Core i7 and 32GB of RAM. For consistency, slack variables are available at all nodes in the system. 
However, slack variables within TPIA can be included at only a subset of system nodes based on real-world system knowledge.
To demonstrate operational bounds, in all test cases we used $0.95$ per unit (p.u.) of the nominal bus voltage for $\overline{V^\phase_i}$ and $1.05$ p.u for $\underline{V^\phase_i}$, a range commonly enforced by real-world utilities.  The diode limiting heuristic is used for all formulations. 

We tested our approaches on six test cases from PNNL taxonomy feeders\cite{PNNL}, which range in size from 76 to 5360 nodes and represent a mix of urban, suburban, and rural systems across varying U.S. geographic regions. These test cases all use constant PQ load models, although equivalent circuit approaches can accommodate other load models. There is a variety of line and three-phase transformer configurations represented. The test cases are available  at \cite{github}. In R1-12.47-3\_OV and R2-25.00-1\_OV, the load factor has been increased compared to the original PNNL case in order to make the cases infeasible such that a standard power flow solver cannot solve them. The other four cases were used as-is, but their power flow solutions exhibit low bus voltages.

\section{Results \& Discussion} \label{result}

In this section, we demonstrate the efficacy of TPIA for solving real-world distribution system problems faced by grid planners. In the first set of results, we show the value of TPIA in determining the location and amount of \textit{extra power} needed to reach a feasible power flow solution within acceptable voltage margins. In the second set of results, we demonstrate how grid planners can use \B{} to ensure network feasibility by deploying only reactive power assets.
\subsection{Enforcing Voltage Bounds on Large-Scale Feeders}
Table \ref{table:pq} presents results from applying \I{} from (\ref{eqn:TPIA_I}) to six test cases of increasing size; for comparison, we also report the results produced by our three-phase power-flow solver described in \cite{Amrit}, which cannot enforce voltage limits since it is only a simulation.
Although our objective function is the sum of squares of complex power slack injections, the values reported in the $\Sigma |P_{feas}|$ and $\Sigma |Q_{feas}|$ columns are the sum of absolute values of power injections, a more useful metric for quantifying the compensation needed for corrective action.

In Table \ref{table:pq}, three-phase power flow did not converge (DNC) for the overloaded cases R1-12.47-3\_OV and R2-25.00-1\_OV. A solution could not be determined in these power deficient networks where the system slack bus is the only generation source.
Power flow does converge for the bottom four test cases of Table \ref{table:pq}, but each solution contains buses with at least one phase with a voltage magnitude ($V_{mag}$) below 95\% of its nominal value.  
The amount of injected power needed in each TPIA solution is correlated with the number of voltage violations in the power flow solution. 
Because TPIA is an optimization and requires Lagrangian multiplier variables, the problem size is much larger than the same network in power flow.
As a result, TPIA requires more iterations and time to converge than power flow. 
However, the TPIA solution times scale linearly with the test case size.

\begin{figure}[tbp]
    \centering
    \includegraphics[width = 2.8in]
    {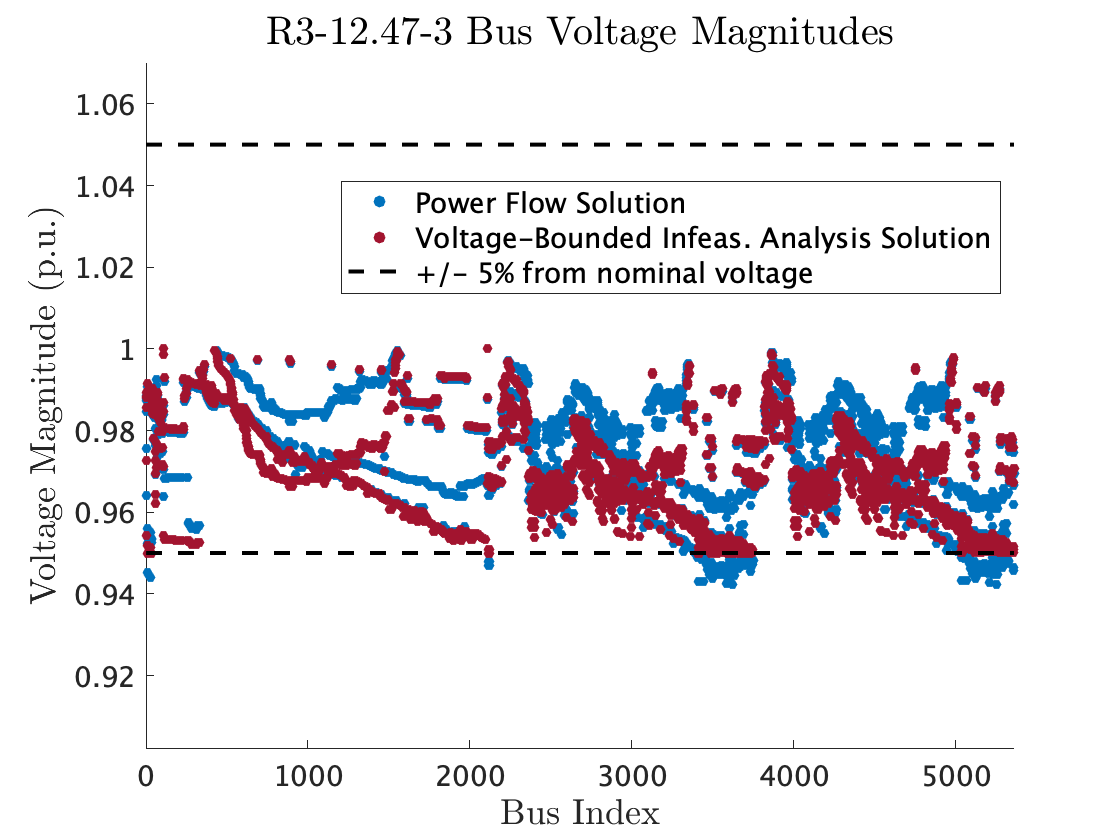}
    \squeezeup \squeezeup 
    \caption{The impact of slack  injections to satisfy $\pm 0.05$ p.u. bus voltage magnitude limits for R3-12.47-3.}
    \label{fig:busvoltage}
    \squeezeup \squeezeup
\end{figure}

To illustrate more explicitly how infeasibility sources and operational bounds can affect the feeder voltages across networks, Fig. \ref{fig:busvoltage} shows the per-phase voltage magnitudes of all nodes in the R3-12.47-3 test case for two solutions. 
A three-phase power flow solution for this case contains 349 bus voltages below 0.95 of the nominal bus voltage (shown by the blue dots), and we observe that the low voltage issues become especially pronounced at the ends of the three radial areas of the network. 
The TPIA solution (shown by the red dots), keeps these low-voltage nodes within the defined bounds through the deployment of infeasibility injections representing 12\% and 11\% of the total P and Q consumed by the feeder, respectively. To implement the output of this analysis, a combination of DERs and reactive power assets could be deployed at the nodes where the slack sources are active in the solution. 

\subsection{Quantifying and Siting Reactive Power Compensation}
Voltage-bounded solutions were also found using \B{}, where only reactive slack sources are available for compensation. 
The results are shown in Table \ref{table:q-only}. 
Across all cases, restricting the slack sources results in at most three more iterations and at most 28\% longer run time. 
Additionally, more total reactive power is injected than the total $Q$ from \I{} in Table \ref{table:q-only}, but now there is zero real power deployed.
\begin{table}[h]
\squeezeup \squeezeup
\centering
\caption{\label{table:q-only} Results when only reactive slack injections are allowed}
\begin{tabular}{|c|c|c|c|c|}
\hline
\textbf{Test case} & \textbf{$\Sigma |P_{feas}|$} & \textbf{$\Sigma |Q_{feas}|$ (kVAR)} & \textbf{Iter.} & \textbf{Time (s)} \\ \hline
R1-12.47-3\_OV     & 0    & 909.90                & 10             & 0.57              \\ \hline
R2-25.00-1\_OV     & 0    & 19458.10              & 27             & 14.33             \\ \hline
R4-12.47-1         & 0    & 1535.38               & 32             & 34.83             \\ \hline
R1-12.47-1         & 0    & 1639.78               & 22             & 26.04             \\ \hline
R5-12.47-3         & 0    & 6356.72               & 32             & 91.29             \\ \hline
R3-12.47-3         & 0    & 1142.25               & 35             & 151.04            \\ \hline
\end{tabular}
\squeezeup
\end{table}

To demonstrate using TPIA for reactive power compensation planning, we selected R1-12.47-3\_OV, which is infeasible when run with a standard three-phase power flow. We ran \B{} to enforce enforce voltage limits and make the network feasible using only reactive power compensation on nodes at or above 480V and found that 0.6 MVAR of reactive power was needed. We then modified the network to add capacitors using the \B{} solution for guidance. The modified case was validated by running it in three-phase power flow, which converged in 4 iterations with all bus voltage magnitudes between 0.95 and 1.05 p.u. 
\section{Conclusion}
Three-phase analysis tools used for distribution grid planning must account for the growing distribution grid complexity. We developed a robust three-phase infeasibility analysis (TPIA) framework with operational bounds that improves existing planning capabilities by providing realistic actionable solutions. We presented TPIA with slack current, reactive power, and susceptance models to develop methods that not only identify the extent of infeasibility in a network, but also identify and localize the real-world assets needed to accomplish grid feasibility. TPIA that uses slack reactive power or admittance provides reactive power compensation and represents specific grid assets like STATCOMs and capacitors. Including voltage bounds pushed solutions towards realistic operating ranges and demonstrated the ability of slack sources to remedy not just situations where power flow fails, but also those where power flow converges to unrealistic solutions. 


\vspace{12pt}


\begin{thebibliography}{00}
\bibitem{IEA} International Energy Agency, ``Global EV Outlook 2021", April 2021. 

\bibitem{NAS} National Academy of Sciences, Engineering, and Medicine. \textit{The Future of Electric Power in the United States}. Washington, D.C: the National Academies Press, 2021. 

\bibitem{derreport} J. S. Homer, et. al. ``Electric Distribution System Planning with DERs- High-level Assessment of Tools and Methods",  Pacific Northwest National Lab, Richland, WA, March 2020.

\bibitem{SCE} ``Reimagining the Grid". Southern California Edison, Rosemead, CA, USA. Dec. 2020 

\bibitem{DGcommittee} A. Keane et al., ``State-of-the-Art Techniques and Challenges Ahead for Distributed Generation Planning and Optimization," \textit{IEEE Trans. on Power Sys.}, vol. 28, no. 2, May 2013.

\bibitem{gill} S. Gill, I. Kockar, and G. Ault, ``Dynamic Optimal Power flow for Active Distribution Networks," \textit{IEEE Trans. on Power Sys.}, vol. 29, no. 1, Jan. 2014.

\bibitem{Dall} A. S. Zamzam, N. D. Sidiropoulos, and E. Dall'Anese, ``Beyond Relaxation and Newton-Raphson: Solving AC-OPF for Multi-Phase Systems with Renewables," \textit{IEEE Trans. on Smart Grid}, vol. 9, no. 5, Sept. 2016.

\bibitem{Ochoa} L. Ochoa and G. Harrison, ``Minimizing Energy Losses: Optimal Accommodation and Smart Operation of Renewable Distributed Generation," \textit{IEEE Trans.on Power Sys.}, vol. 26, no. 1, Feb. 2011.

\bibitem{Claeys} S. Claeys, F. Geth, and G. Deconinick, ``Optimal Power Flow in Four-Wire Distribution Networks: Formulation and Benchmarking," \textit{Electric Power Sys. Res.}, vol. 213, Dec. 2022.

\bibitem{PSCC} E.Foster, A. Pandey, and L. Pileggi, ``Three-Phase Infeasibility Analysis for Distribution Grid Studies," \textit{Electric Power Sys. Res.}, vol. 212, Nov. 2022.

\bibitem{overbye} T. J. Overbye,``A power flow measure for unsolvable cases,'' \textit{IEEE Trans. on Power Systems}, vol. 9, no. 3, Aug. 1994.

\bibitem{marko} M. Jereminov, D. M. Bromberg, A. Pandey, M. R. Wagner and L. Pileggi, ``Evaluating Feasibility Within Power Flow," \textit{IEEE Trans.on Smart Grid}, vol. 11, no. 4, July 2020.

\bibitem{Cindy}  S. Li, A. Pandey, A. Agarwal, M. Jereminov and L. Pileggi, "A LASSO-Inspired Approach for Localizing Power System Infeasibility," presented at \textit{2020 IEEE PES GM}, Montreal, QC, Canada, Aug. 2-6, 2020.

\bibitem{homotopy} A. Pandey, M. Jereminov, M. R. Wagner, G. Hug and L. Pileggi, "Robust Convergence of Power Flow Using TX Stepping Method with Equivalent Circuit Formulation," presented at the Power Sys. Comp. Conf. (PSCC), Dublin, Ireland, June 11-15, 2018.

\bibitem{diode} M. Jerminov, A. Pandey, and L. Pileggi, ``Equivalent Circuit Formulation for Solving AC Optimal Power Flow," \textit{IEEE Trans. on Power Sys.}, vol. 34, no. 3, May 2019.

\bibitem{shunt_max} \textit{IEEE Standard for Shunt Power Capacitors}, IEEE Std 18-2012 (Revision of IEEE Std 18-2002), Feb. 15, 2013.

\bibitem{Amrit}  A. Pandey, et al., ``Robust Power Flow and Three-Phase Power Flow Analyses," in \textit{IEEE Trans. on Power Sys.}, vol. 34, no. 1, Jan. 2019.

\bibitem{Boyd} S. Boyd, S. P. Boyd, and L. Vandenberge, \textit{Convex Optimization}, Cambridge University Press, 2004.

\bibitem{two_stage} T. McNamara, A. Pandey, A. Agarwal, and L. Pileggi, ``Two-stage homotopy method to incorporate discrete control variables into AC-OPF," \textit{Electric Power Sys. Res.}, vol. 212, Nov. 2022.

\bibitem{code} A. Pandey, N. Turner-Bandele, E. Foster, and T. McNamara, \textit{SUGAR-D: A Distribution Systems Analysis and Optimization Tool v2.02}, CMU ECE Pileggi Lab Group, Pittsburgh, PA, 2022. 

\bibitem{PNNL} K.P. Schneider, et al., ``Modern Grid Initiative Distribution Taxonomy Final Report," Pacific Northwest National Lab, Nov. 2008.

\bibitem{github} E. Foster, TPIA Experiment Test Cases, 2022, Github Repository. https://github.com/emfoster/TPIA-Experiment-Testcases


\end{thebibliography}
\end{document}